\documentclass[12pt]{amsart}
\usepackage{amsmath}
\usepackage{amsfonts,amssymb,amsthm}
\newtheorem{theorem}{Theorem}
\newtheorem{lemma}[theorem]{Lemma}
\newtheorem{corollary}[theorem]{Corollary}
\begin{document}
\title[Alternating groups]{A property of alternating groups}
\author{Henry Cejtin and Igor Rivin}
\address{Sourcelight Technologies, Evanston, IL}
\email{henry@sourcelight.com}
\address{Department of Mathematics, Temple University, Philadelphia}
\curraddr{Mathematics Department, Princeton University}
\email{rivin@math.temple.edu}
\date{today}
\begin{abstract}
We describe an efficient algorithm to write any element of the
alternating group $A_n$ as a product of two $n$-cycles (in particular,
we give a simple proof that any element of $A_n$ can be so written). An easy
corollary is that every element of $A_n$ is a commutator in $S_n$ 
\end{abstract}
\maketitle

Consider the alternating group $A_n.$ We will show the following 
\begin{theorem}
\label{main}
There is an $O(n)$ algorithm to write any element $\sigma \in A_n$ 
as a product of two $n$-cycles.
\end{theorem}

That any element of $A_n$ \textit{can} be so written is a result of
E.~A.~Bertram \cite{bertie}, but his argument is less direct and hence
does not seem to lead to an optimal algorithm. On the other hand, Bertram's
result is more general. 

A corollary (Corr.~\ref{commu}), it is seen that every
element of $A_n$ is a commutator in $S_n$ -- the well-known
fact that $A_n$ is the commutator subgroup of $S_n$ gives a much
weaker statement. As an immediate corollary of the (one line) proof of
Corr.~\ref{commu} we obtain an $O(n)$ algorithm to represent any
element of $A_n$ as a commutator.

In the sequel we use exclusively the cycle notation for
permutations. Both the existence and the algorithmic aspects of the
cycle decomposition are discussed in \cite[Section 1.3.3]{knuthv1}.

In order to prove Theorem \ref{main} we will need a few lemmas:

\begin{lemma}
\label{intertwine}
Let $\sigma_1 = (a_1, \dots, a_s),$ and $\sigma_2 = (b_1, \dots,
b_t).$ Furthermore suppose that $\{a_1, \dots, a_s\} \cap \{b_1, \dots, b_t\} =
\emptyset.$ 
Let $\rho = (a_s, b_t).$ Then $\tau_1 = \sigma_1 \sigma_2 \rho$ is an $s+t$-cycle,
as is $\tau_2 = \rho \sigma_1 \sigma_2.$
\end{lemma}

\begin{proof}
Since $\tau_2 = \rho \tau_1 \rho^{-1},$ it is enough to show that
$\tau_1$ is an $s+t$-cycle. Since the cycle decomposition is conjugacy
invariant, we can assume without loss of generality that 
\begin{gather*}
\sigma_1 = (1, 2, \dots, s),\\
\sigma_2 = (s+1, \dots, s+t),\\
\rho = (s, s+t).
\end{gather*}
Then
\begin{equation*}
\tau_1(x) = \begin{cases}
x+1, \qquad x < s-1,\\
s+t, \qquad x = s-1,\\
s+1, \qquad x = s+t,\\
x+1, \qquad s < x < s+t -1,\\
s, \qquad, x = s+t - 1,\\.
x, \qquad x > s+t
\end{cases}
\end{equation*}
The assertion of the lemma follows by inspection.
\end{proof}

\begin{lemma}
\label{induction}
Suppose $\sigma \in S_n,$ where $n=t+s,$ can be written as $\sigma_1
\sigma_2,$ where 
$\sigma_1$ only acts non-trivially on $\{1, \dots, t\}$ and $\sigma_2$
only acts non-triviall on $\{t+1, \dots, t+s\}.$
Suppose further that $\sigma_1$ is a product of two
$t$-cycles and $\sigma_2$ is a product of two $s$-cycles. Then
$\sigma$ is a product of two $n$-cycles.
\end{lemma}

\begin{proof}
Let $\tau = (t, t+s).$ Let $\sigma_1 = \rho_{11} \rho_{12},$ while
$\sigma_2 = \rho_{21}\rho_{22},$ where the $\rho_{ij}$ are the $t$- and
$s$-cycles as per hypothesis of the Lemma. Then:
\begin{equation*}
\sigma_1 \sigma_2 = \rho_{11} \rho_{12} \rho_{21} \rho_{22} =
\rho_{11} \rho_{21} \rho_{12} \rho_{22} = 
(\rho_{11} \rho_{21} \tau) (\tau \rho_{12} \rho_{22}).
\end{equation*}
By Lemma \ref{intertwine} both the parenthesized terms are $n$-cycles,
and the result follows.
\end{proof}

\begin{lemma}
\label{onecyc}
Let $n$ be odd. Then 
any $n$-cycle $\rho$ is a product of two $n$-cycles.
\end{lemma}
\begin{proof}
The cyclic subgroup $C_{\rho}$ generated by $\rho$ is of order $n,$ 
$m = \frac{n+1}2$ is an integer, and clearly $\left(\rho^m\right)^2 =
\rho.$ Since $n+1$ is relatively prime to $n,$ $\rho^m$ is a generator
of $C_{\rho},$ hence an $n$-cycle.
\end{proof}

\begin{lemma}
\label{eqlem}
Let $n = 4 m,$ and let $\sigma \in S_n$ be the product of $\rho_1$ and
$\rho_2,$ where $\rho_1$ and $\rho_2$ are $2m$-cycles, and $\rho_1$ and
$\rho_2$ are disjoint. Then $\sigma$ is a product of two $n$-cycles.
\end{lemma}
\begin{proof}
Without loss of generality, we can assume that 
\begin{gather*}
\rho_1 = (1, 3, \dots, 4m - 1),\\
\rho_2 = (2, 4, \dots, 4m).
\end{gather*}
It is easy to see that 
$\sigma = \rho^2,$ where 
\begin{equation*}
\rho = (1, 2, \dots, 4m)
\end{equation*}
\end{proof}

\begin{lemma}
\label{uneq}
Let $n = 2s + 2t,$ where $s < t.$ Let $\sigma \in S_n$ be the
product of $\rho_1$ and $\rho_2,$ where $\rho_1$ is a $2s$-cycle and
$\rho_2$ is a $2t$-cycle, and $\rho_1$ is disjoint from $\rho_2.$ Then
$\sigma$ is a product of two 
$n$-cycles.
\end{lemma}
\begin{proof}
First, without loss of generality, assume that
\begin{gather*}
\rho_1 = (1, 3, \dots, 4s -1),\\
\rho_2 = (2, 4, 6, \dots, 2s + 2t, 4s +1, 4s+3, \dots,
2s + 2 t -1).
\end{gather*}

Now, define $\lambda_1, \lambda_2 \in S_n$ as follows:
\begin{equation*}
\lambda_1(x) = \begin{cases}
4s+1 \qquad x=1,\\
1    \qquad x = 4s,\\
2    \qquad x = 2s + 2 t,\\
x+1 \qquad \mbox{otherwise}.
\end{cases}
\end{equation*}
In other words, $\lambda_1$ is the $n$-cycle
\[(1, 4s + 1, 4s + 2, \dots, 2s + 2t, 2, 3, \dots, 4s).\]
\begin{equation*}
\lambda_2(x) = \begin{cases}
1 \qquad x = n,\\
x+1 \qquad \mbox{otherwise}.
\end{cases}
\end{equation*}
In other words, $\lambda_2$ is the $n$-cycle
\[(1, 2, 3, \dots, n).\]
Now, let us compute $\lambda_1 \lambda_2:$
\[\lambda_1 \lambda_2(n) = \lambda_1(1) = 4s + 1.\]
If $x \neq n,$
\[\lambda_1 \lambda_2(x) = \lambda_1(x+1) = \begin{cases}
1 \qquad x = 4s -1,\\
2 \qquad = n -1,\\
x+2 \qquad \mbox{otherwise}.
\end{cases}\]
It follows immediately that 
\[\lambda_1 \lambda_2 = \rho_1 \rho_2.\]
\end{proof}
\begin{proof}[Proof of Theorem \ref{main}]
The proof will proceed by induction on $n.$ If $n=1,$ the
result is obvious. Now, consider $\sigma \in A_n.$ If $\sigma =
\sigma_1 \sigma_2,$ where $\sigma_1$ and $\sigma_2$ are even and
disjoint, the assertion of Theorem \ref{main} will follow from the
induction hypothesis and Lemma \ref{induction}. 

Consider now the cycle decomposition of $\sigma$ -- when we count the
cycles below, we \emph{do} count trivial (one-element) cycles.
\begin{itemize}
\item{\textbf{The cycle decomposition of $\sigma$ has exactly one cycle.}}  In
this case, the statement of the Theorem follows from Lemma
\ref{onecyc}.

\item{\textbf{The cycle decomposition of $\sigma$ has exactly two cycles.}}
The statement of the Theorem follows from Lemmas \ref{eqlem} and
\ref{uneq}.

\item{\textbf{The cycle decomposition of $\sigma$ has more than two cycles.}}
If one of the cycles is of odd length, then it is an even
permutation, and the result follows from the discussion at the
beginning of the proof. If no cycle is of odd length, then, since
$\sigma$ is an even permutation, there is an even number of cycles in
the cycle decomposition of $\sigma,$ so
\[\sigma = \rho_1 \rho_2 \dots \rho_{2k} = (\rho_1 \rho_2) (\rho_3
\dots \rho_{2k}).\]
Each of the parenthesized terms is an even permutation, and so
$\sigma$ is a product of two $n$-cycles by the discussion at the
beginning of this proof.
\end{itemize}
\end{proof}

\begin{corollary}
\label{commu}
Every element of $A_n$ is a commutator of a pair of elements in $S_n.$
\end{corollary}

\begin{proof}
All $n$-cycles are conjugate, and if $\rho$ is an $n$-cycle, so is
$\rho^{-1}.$ Now let $\sigma \in A_n.$ By Theorem \ref{main}, $\sigma
= \rho_1 \rho_2,$ where $\rho_1, \rho_2$ are $n$-cycles. Since $\rho_2
= \tau \rho_2^{-1} \tau^{-1},$ we see that 
$\sigma = [\rho_1, \tau].$
\end{proof}
\bibliographystyle{alpha}

\end{document}